\documentclass{article}
\usepackage[margin=1.5in]{geometry}
\usepackage{graphicx} 
\usepackage{amsthm}
\usepackage{amssymb}
\usepackage{amsmath} 
\usepackage{mathtools}
\usepackage{xcolor}
\usepackage{enumerate}
\usepackage{microtype}
\usepackage{hyperref}
\usepackage[natbib,style=numeric,sorting=nyt]{biblatex}
\newtheorem{theorem}{Theorem}

\newtheorem{lemma}[theorem]{Lemma}

\newtheorem{cor}[theorem]{Corollary}

\newtheorem{example}{Example}

\if10     
\usepackage[mathlines]{lineno}
\newcommand*\patchAmsMathEnvironmentForLineno[1]{%
  \expandafter\let\csname old#1\expandafter\endcsname\csname #1\endcsname
  \expandafter\let\csname oldend#1\expandafter\endcsname\csname end#1\endcsname
  \renewenvironment{#1}%
     {\linenomath\csname old#1\endcsname}%
     {\csname oldend#1\endcsname\endlinenomath}}%
\newcommand*\patchBothAmsMathEnvironmentsForLineno[1]{%
  \patchAmsMathEnvironmentForLineno{#1}%
  \patchAmsMathEnvironmentForLineno{#1*}}%
\AtBeginDocument{%
\patchBothAmsMathEnvironmentsForLineno{equation}%
\patchBothAmsMathEnvironmentsForLineno{align}%
\patchBothAmsMathEnvironmentsForLineno{flalign}%
\patchBothAmsMathEnvironmentsForLineno{alignat}%
\patchBothAmsMathEnvironmentsForLineno{gather}%
\patchBothAmsMathEnvironmentsForLineno{multline}%
}
\linenumbers
\fi

\def\inst#1{$^{#1}$}

\title{Proof of the Harary--Hill conjecture: an exposition\thanks{The author was supported by grant no.\ 23-04949X of the Czech Science Foundation (GA\v{C}R).
}}
\author{
Martin Balko\inst{1}
}

\begin{document}

\maketitle

\begin{center}
{\footnotesize
\inst{1} 
Department of Applied Mathematics, \\
Faculty of Mathematics and Physics, Charles University, Czech Republic \\
\texttt{balko@kam.mff.cuni.cz}
}
\end{center}

\begin{abstract}
The Harary--Hill conjecture states that the crossing number of the complete graph $K_n$ is exactly
$\frac{1}{4}\left\lfloor\frac{n}{2}\right\rfloor
\left\lfloor\frac{n-1}{2}\right\rfloor
\left\lfloor\frac{n-2}{2}\right\rfloor
\left\lfloor\frac{n-3}{2}\right\rfloor$.
We present a streamlined variant of the proof of this conjecture recently announced by OpenAI (2026). 
Although the original proof is relatively short, it is quite difficult to follow. 
We simplify several steps and provide a concise and accessible exposition of the proof. 
We also combine this new result with known reductions to determine the crossing numbers of several other graph families.
\end{abstract}

\section{Introduction}

In a \emph{drawing} $D$ of a graph $G$, vertices of $G$ are represented by distinct points of the plane and edges of $G$ are simple arcs joining their endpoints and avoiding all other vertices.
As usual, we identify the vertices and their images, as well as the edges and the arcs representing them. We require that the edges pass through no vertices other than their endpoints. We also assume for simplicity that any two edges have only finitely many points in common, no two edges touch at an interior point, and no three edges meet at a common interior point.

A \emph{crossing} is a point that lies in the interior of at least two edges.
We write ${\rm cr}(D)$ for the number of crossings in a drawing $D$.
For a graph $G$ the \emph{crossing number of $G$}, denoted by ${\rm cr}(G)$, is the minimum value of ${\rm cr}(D)$ taken over all drawings $D$ of $G$.

According to one of the most famous problems in topological graph theory, called the \emph{Harary--Hill conjecture} (also known as \emph{Guy’s conjecture}), the crossing number of the complete graph $K_n$ on $n$ vertices satisfies
\[{\rm cr}(K_n) = Z(n)\]
for every positive integer $n$, where
\[
Z(n)=\frac{1}{4}\left\lfloor\frac{n}{2}\right\rfloor
\left\lfloor\frac{n-1}{2}\right\rfloor
\left\lfloor\frac{n-2}{2}\right\rfloor
\left\lfloor\frac{n-3}{2}\right\rfloor.
\]
This conjecture dates back to the late 1950s, when the British
artist Anthony Hill constructed drawings of $K_n$ with exactly
$Z(n)$ crossings; see~\cite{BeinekeWilson2010}. 
The conjectured optimality of these drawings
was first published by Guy~\cite{Guy1960} and later discussed
by Harary and Hill~\cite{HararyHill1963}, while the problem was
further popularized by Erd\H{o}s and Guy~\cite{ErdosGuy1973}.

Although the conjecture resisted proof for more than six decades, considerable progress was made for special classes of drawings.
The conjecture was first proved for \emph{2-page drawings}, in which
the vertices lie on a line, and each edge is drawn in one of the
two halfplanes~\cite{AbregoEtAl2013}. Subsequently,
\'Abrego et al.~\cite{AbregoEtAl2014} introduced the more
general class of \emph{shellable drawings} and established the
conjectured lower bound for this class, which includes both cylindrical and $x$-bounded drawings.
Independently, Balko, Fulek, and Kyn\v{c}l~\cite{BalkoFulekKyncl2015} proved the conjecture for
\emph{$x$-monotone drawings}, even for more general notions of the crossing number.
The classes of drawings for which the Harary--Hill conjecture was known to be true were then further extended to \emph{bishellable
drawings}~\cite{AbregoEtAl2018}, and \emph{seq-shellable
drawings}~\cite{MutzelOettershagen2018}.
More recently, Streltsova and Wagner~\cite{StreltsovaWagner2025} established the conjectured
lower bound for \emph{spherical geodesic drawings}, in which
vertices lie on a sphere, and edges follow unique shortest geodesic arcs.

For arbitrary drawings, the conjecture remained open. 
By 2021, it had been verified for
$n\leq14$ using
computer-assisted methods~\cite{Aichholzer2021}.
Significant progress was also made on asymptotic lower bounds.
Using flag algebras, Balogh, Lidick\'y, and
Salazar~\cite{BaloghLidickySalazar2019} proved that
${\rm cr}(K_n)\geq(0.985+o(1))Z(n)$.
This was later improved to $(0.9899+o(1))Z(n)$ by
Balogh et al.~\cite{BaloghEtAl2026}.

\section{Our contributions}

In October 2026, OpenAI announced a proof of the Harary--Hill conjecture
for all $n$~\cite{openAI}, resolving the problem after more than 60 years.
In this note, we present a simplified variant of this proof. 

\begin{theorem}[\cite{openAI}]
\label{thm-HararyHill}
For every positive integer $n$, the crossing number of $K_n$ equals $Z(n)$.
\end{theorem}

The proof is based on the polynomial method and is given in Section~\ref{sec-proof}.
We simplify some steps of the original argument and give a succinct proof, which requires only elementary knowledge of linear algebra and graph theory.
We believe that a simpler and more accessible proof is valuable to the community, as the underlying method may find applications to other problems concerning crossing numbers.
We also note that an independent exposition of OpenAI's proof is being prepared by Gelasio Salazar~\cite{gelasio}, who has given a detailed presentation of the proof, exploring its broader context, the motivation behind its key arguments, and the origins of the underlying ideas.

We also show that the proof of Theorem~\ref{thm-HararyHill} extends to other notions of the crossing number discussed, for example, in~\cite{FulekPelsmajerSchaeferStefankovic2012}.
A drawing $D$ of a graph $G$ is called \emph{semisimple} if
any two adjacent edges intersect only at their common endpoint.
For two distinct edges $e,f\in E(G)$, let $c_D(e,f)$ denote
the number of crossings between $e$ and $f$ in $D$.
The \emph{pair-crossing number with Rule $+$}, denoted by
${\rm pcr}_+(G)$, is the minimum number of unordered pairs
of crossing edges over all semisimple drawings of $G$.
Similarly, the \emph{odd-crossing number with Rule $+$},
denoted by ${\rm ocr}_+(G)$, is the minimum number of unordered
pairs of edges crossing an odd number of times over all
semisimple drawings of $G$.
Equivalently,
\[
{\rm pcr}_+(G)
=\min_D
\left|\left\{\{e,f\}\in\binom{E(G)}{2}:c_D(e,f)>0\right\}\right|
\]
and
\[
{\rm ocr}_+(G)
=\min_D
\left|\left\{\{e,f\}\in\binom{E(G)}{2}:
c_D(e,f)\equiv1\pmod 2\right\}\right|,
\]
where both minima range over all semisimple drawings $D$ of $G$.
In Subsection~\ref{subsec-variants}, we show that the proof actually gives
\[{\rm ocr}_+(K_n) = {\rm pcr}_+(K_n) = Z(n)\]
for every positive integer $n$.

As immediate consequences, we also include a list of several graph classes whose crossing numbers are determined exactly using Theorem~\ref{thm-HararyHill} and known reductions.
This includes complete graphs with a deleted perfect matching or a single deleted edge, Cartesian products of complete graphs with paths or even cycles, and extended Sierpi\'nski graphs; see Section~\ref{sec-applications}.

\section{Proof of Theorem~\ref{thm-HararyHill}}
\label{sec-proof}

In this section, we give the proof of the Harary--Hill conjecture.
That is, we prove Theorem~\ref{thm-HararyHill}.
We start by giving a brief overview of the whole proof.

\paragraph{Proof overview.} 
The proof combines Tutte's algebraic approach to crossing numbers~\cite{Tutte1970}
with a basic principle of the polynomial method: if a finite set $\mathcal{X}$ has
fewer elements than the dimension of a vector space $\mathcal S$ of polynomials,
then there exists a nonzero polynomial in $\mathcal S$ that vanishes at every point of $\mathcal{X}$.
Thus, to obtain the lower bound $|\mathcal{X}|\ge\dim\mathcal S$, it suffices to show that the only polynomial in $\mathcal S$ vanishing at every point of $\mathcal{X}$ is the zero polynomial.
This dimension-counting principle is a basic ingredient of the polynomial method; see, for example,~\cite{Dvir2009}. 
Related polynomial interpolation
and nonvanishing techniques also play a central role in recent breakthroughs on Landau--Siegel zeros~\cite{Green2026Siegel,OpenAI2026Siegel}.
For us, $\mathcal{X}$ will be a set of pairs of crossing edges of a simple drawing $D$ of $K_n$ with $n=2s+1$, and $\mathcal{S}$ will be the vector space of polynomials $P(u,v,x,y)$ of degree at most $s-2$ separately in each variable that are symmetric within $(u,v)$ and $(x,y)$.
Then, $\dim \mathcal{S} = \binom{s}{2}^2=Z(2s+1)$.
The connection between $\mathcal{X}$ and $\mathcal{S}$ comes from Tutte's approach.
We consider edges $[v_iv_j]$ of an oriented $\vec{K}_n$ as elements of a basis of a real vector space $C_1(\vec{K}_n;\mathbb{R})$.
\emph{Cycles} are then elements of $C_1(\vec{K}_n;\mathbb{R})$ whose incoming and outgoing coefficients balance at each vertex.
Representing each vertex $v_i$ with a unique real number $t_i$, we use Lagrange interpolation to obtain nonzero weights $w_i$ such that expressions
$\sum_{i,j}w_iw_jt_i^at_j^b[v_iv_j]$
are cycles for all $0\le a,b\le n-2$.
Let $\mathcal{X}_{\pm}$ be the set of all ordered pairs of undirected edges that cross in $D$, and note that $|\mathcal{X}_{\pm}| = 2{\rm cr(D)}$.
Using the result about cycles, we show that there are nonzero coefficients $a_\alpha$ such that every polynomial $P \in \mathcal{S}$ satisfies $\sum_{\alpha \in \mathcal{X}_{\pm}} a_\alpha P^2(t_i,t_j,t_k,t_l) = 0$ where $\alpha = (\{v_i,v_j\},\{v_k,v_l\})$.
This, together with $a_{(e,f)} = -a_{(f,e)}$, implies that the set $\mathcal{X} = \{\alpha \in \mathcal{X}_{\pm} \colon a_\alpha > 0\}$ has size ${\rm cr}(D)$.
Thus, if ${\rm cr}(D)$ is smaller than $Z(2s+1) = \dim \mathcal{S}$, then $|\mathcal{X}| < \dim \mathcal{S}$, and the homogeneous system  $P(t_i,t_j,t_k,t_l) =0$ for every $(\{v_i,v_j\},\{v_k,v_l\}) \in \mathcal{X}$ has a nonzero solution $P\in\mathcal S$.
We choose one such $P$ of minimum degree in the first variable.
However, a second application of the identity for cycles,
together with an interpolation argument, shows that
$P(u,v,u,y)\equiv0$.
By the symmetries of $P$, this implies that $C=(u-x)(u-y)(v-x)(v-y)$ divides $P$.
Since $C$ is nonzero at all evaluations corresponding to pairs of disjoint edges, $P/C$ is a nonzero polynomial
in $\mathcal S$ that still vanishes on $\mathcal X$ but has smaller degree in $u$ than $P$, contradicting the choice of $P$.
Thus, ${\rm cr}(D) \geq Z(2s+1)$.
The cases for even $n$ follow by standard double counting from the odd cases.

\subsection{Algebraic preliminaries}

Here, we state a few interpolation results and definitions that are needed later in the proof of Theorem~\ref{thm-HararyHill}.
All results in this subsection are purely algebraic and do not involve graph drawings.

The following consequence of Lagrange interpolation assigns nonzero weights $w_i$ to $n$ distinct real numbers $t_1,\dots,t_n$ so that every polynomial $P(x)$ of degree at most $n-2$ has zero weighted sum $\sum_{i=1}^n w_i P(t_i)$, a cancellation property we will later use.

\begin{lemma}
\label{lem-3.1}
For a positive integer $n$ and distinct real numbers $t_1,\dots,t_n$, let $F(x) = \prod_{i=1}^n (x-t_i)$.
Then, every real polynomial $P(x)$ of degree at most $n-2$ satisfies
\[
\sum_{i=1}^n \frac{P(t_i)}{F'(t_i)} = 0.
\]
\end{lemma}
\begin{proof}
Using Lagrange interpolation~\cite{StoerBulirsch2002}, we can express $P$ in terms of its values at $t_1,\dots,t_n$ and obtain
\[
P(x) = \sum_{i=1}^n P(t_i) \prod_{j \neq i} \frac{x-t_j}{t_i-t_j} =\sum_{i=1}^n
\frac{P(t_i)}{F'(t_i)}
\prod_{j\neq i}(x-t_j),
\]
where the second equality follows from $F'(t_i)=\prod_{j\neq i}(t_i-t_j)$.
The coefficient of $x^{n-1}$ is zero in $P(x)$ as $\deg P\le n-2$.
On the other hand, each product on the right side is monic of degree $n-1$ and the coefficient of $x^{n-1}$ is $\sum_{i=1}^n
\frac{P(t_i)}{F'(t_i)}$.
Thus, we obtain
$
\sum_{i=1}^n\frac{P(t_i)}{F'(t_i)}=0.
$
\end{proof}

For non-negative integers $a$ and $b$, we let $\mathcal{P}_{a,b}$ denote the vector space of real polynomials in $x$ and $y$ of degrees at most $a$ in $x$ and $b$ in $y$, respectively.
The following lemma shows that arbitrary values on two triangular sets containing $s(s-1)$ grid points in total uniquely specify a polynomial in either of two prescribed spaces $\mathcal{P}_{a,b}$, where $(a,b)\in\{(s-1,s-2),(s-2,s-1)\}$.

\begin{lemma}
\label{lem-3.2}
For an integer $s\ge2$, let $t_1,\dots,t_{2s}$ be distinct real numbers.
Set
\[
L_s=
\{(i,j):1\le j<i\le s\}
\cup
\{(i,j):s+1\le i<j\le 2s\}.
\]
For either choice
$(a,b)\in\{(s-1,s-2),(s-2,s-1)\}$
and any prescribed real numbers $c_{i,j}$ with $(i,j)\in L_s$,
there exists a unique polynomial $P\in\mathcal P_{a,b}$ satisfying
$P(t_i,t_j)=c_{i,j}$
for all $(i,j)\in L_s$.
In particular, every nonzero polynomial in either space is nonzero
at some point $(t_i,t_j)$ with $(i,j)\in L_s$.
\end{lemma}

\begin{proof}
Consider the system of $s(s-1)$ linear
equations $P(t_i,t_j)=c_{i,j}$ where the variables are the $(a+1)(b+1)=s(s-1)$ coefficients of $P$.
It suffices to show that a polynomial in either space $\mathcal{P}_{a,b}$
vanishing at all points $(t_i,t_j)$ with $(i,j)\in L_s$ is identically zero, as having only
the zero solution makes the coefficient
matrix invertible and guarantees a unique solution for any prescribed values $c_{i,j}$.

We prove this simultaneously for both spaces by induction on $s$.
Interchanging the variables and swapping the two blocks
$t_1,\dots,t_s$ and $t_{s+1},\dots,t_{2s}$ interchanges the two cases.
Thus, it suffices to consider $P\in\mathcal P_{s-1,s-2}$.

For $s=2$, the polynomial $P$ depends only on $x$, has degree
at most one, and vanishes at $t_2$ and $t_3$, so $P=0$.
Now suppose $s\ge3$.
Each of the polynomials $P(t_s,y)$ and $P(t_{s+1},y)$
has degree at most $s-2$ and vanishes at $s-1$ distinct numbers $t_i$.
Hence both are identically zero, and
$
P(x,y)=(x-t_s)(x-t_{s+1})R(x,y)
$
for some $R\in\mathcal P_{s-3,s-2}$.
Delete the numbers $t_s$ and $t_{s+1}$ and relabel the remaining
numbers $t_i$ in their original order.
The remaining pairs in $L_s$ form $L_{s-1}$, and the two
removed factors are nonzero at every corresponding point.
Thus, $R$ vanishes at all these points.
Since
$
R\in\mathcal P_{s-3,s-2}
=\mathcal P_{(s-1)-2,(s-1)-1},
$
the induction hypothesis for $s-1$ gives $R=0$, and consequently $P=0$.
\end{proof}

Let $s\ge2$, let $t_1,\dots,t_{2s}$ be distinct real numbers, and set $F(x)=\prod_{i=1}^{2s}(x-t_i)$.
Let $\prec$ be a linear ordering on $\{t_1,\dots,t_{2s}\}$.
We define
\[
S(i,j)=
\begin{cases}
  \frac{1}{2},  & t_i \prec t_j,\\
  -\frac{1}{2}, & t_j \prec t_i,\\
  0,  & i=j,
\end{cases}
\]
for all integers $i$ and $j$ with $1 \leq i,j \leq 2s$.
For any two real polynomials $P(x,y)$ and $Q(x,y)$, we also set
\begin{equation}
\label{eq-pairing}
M(P,Q) = \sum_{i=1}^{2s}\sum_{j=1}^{2s} S(i,j) \frac{t_i-t_j}{F'(t_i)F'(t_j)} P(t_i,t_j)Q(t_i,t_j).
\end{equation}

The last result in this subsection shows that $M$ detects every
nonzero polynomial $P$ in the prescribed space through a suitable
test polynomial $Q$, regardless of the chosen order $\prec$ of the numbers $t_1,\dots,t_{2s}$.

\begin{lemma}
\label{lem-3.3}
Let $s\ge 2$ be an integer.
If $P\in\mathcal P_{s-2,s-2}$ satisfies
$M(P,Q)=0$ for every $Q\in\mathcal P_{s-1,s-1}$,
then $P=0$.
\end{lemma}

\begin{proof}
We relabel the numbers $t_1,\dots,t_{2s}$ in their prescribed
order $\prec$, so that $S(i,j)=1/2$ for $i<j$ and
$S(i,j)=-1/2$ for $i>j$.
Let $P\in\mathcal P_{s-2,s-2}$ satisfy the hypothesis, and take
$Q\in\mathcal P_{s-2,s-1}\subseteq\mathcal P_{s-1,s-1}$.

For each fixed $i$, the polynomial
$(t_i-y)P(t_i,y)Q(t_i,y)$ has degree at most
$1+(s-2)+(s-1)=2s-2$.
Applying Lemma~\ref{lem-3.1} with $n=2s$ and multiplying
by $1/F'(t_i)$, we obtain
\[
\sum_{j=1}^{2s}
\frac{t_i-t_j}{F'(t_i)F'(t_j)}
P(t_i,t_j)Q(t_i,t_j)=0.
\]
Consequently, subtracting $1/2$ from every $S(i,j)$ with
$i\le s$, and adding $1/2$ to every $S(i,j)$ with $i>s$,
does not change $M(P,Q)$.
Indeed, each contribution for the $i$th term changes by a constant
multiple of the zero sum above.
Since the terms with $i=j$ vanish, this yields
\[
M(P,Q)
={}-\sum_{1\le j<i\le s}
\frac{t_i-t_j}{F'(t_i)F'(t_j)}
P(t_i,t_j)Q(t_i,t_j)
+\sum_{s<i<j\le 2s}
\frac{t_i-t_j}{F'(t_i)F'(t_j)}
P(t_i,t_j)Q(t_i,t_j).
\]
Every coefficient $\frac{t_i-t_j}{F'(t_i)F'(t_j)}$ in these sums is nonzero, since
$i\neq j$ and the nodes are pairwise distinct.

Suppose for contradiction that $P\neq0$.
Then, since $P\in\mathcal P_{s-1,s-2}$,
Lemma~\ref{lem-3.2} gives $(i_0,j_0)\in L_s$ such that $P(t_{i_0},t_{j_0})\neq0$.
The same lemma provides $Q\in\mathcal P_{s-2,s-1}$ whose value $Q(t_i,t_j)$
is one at $(t_i,t_j) = (t_{i_0},t_{j_0})$ and zero at all other points $(t_i,t_j)$ indexed by $L_s$.
The resulting sum has exactly one nonzero term, contradicting
$M(P,Q)=0$. 
Hence $P=0$.
\end{proof}

\subsection{Treating crossings algebraically}

A drawing $D$ of a graph $G$ is \emph{simple} if any two edges intersect at most once.
That is, any two edges can share only either a crossing or a common endpoint.
We say that a drawing $D$ is \emph{polygonal} if all edges are represented by polygonal curves that are straight near every vertex and crossing.
It is a well-known fact that every drawing of $K_n$ with the minimum number of crossings can be redrawn without increasing the number of crossings so that the resulting drawing is simple~\cite{Guy1972}.
By local modifications near vertices and crossings and polygonal
approximation of the remaining arcs, it can also be made polygonal without changing the number of crossings.

We use an approach similar to Tutte's algebraic approach to crossing numbers~\cite{Tutte1970}.
For a polygonal drawing $D$ of an oriented graph, we define a bilinear form $J_D$ that combines signed crossings with local contributions determined by the cyclic orders of edges
at vertices to express the classical fact that two closed oriented curves in the plane have
algebraic intersection number zero; see, for example, \cite[Section~1.3]{Skopenkov2020}.
A closely related mapping appears in~\cite[Lemma~8.1]{GoncharovKenyon2013}.

Let $D$ be a polygonal drawing of an oriented graph $\vec{G}=(V,E)$.
The \emph{real edge space} $C_1(\vec{G};\mathbb R)$ is the
vector space having these oriented edges as a basis. 
Its elements, called \emph{chains}, are formal linear combinations $z=\sum_{e\in E}z_e e$ with real coefficients, where each edge $e$ is taken with its reference orientation. 
For adjacent vertices $u,v$, we write
$[uv]$ for the edge directed from $u$ to $v$, with the convention
$[uv]=-[vu]$.
We also set $[vv]=0$. 
Thus, reversing the direction of an
edge negates its vector. 
The \emph{vertex space} $C_0(\vec{G};\mathbb R)$ has the formal vertex symbols $\{e_v:v\in V\}$ as a basis.
The \emph{boundary map} $\partial:C_1(\vec{G};\mathbb R)\to C_0(\vec{G};\mathbb R)$ is the linear map given by $\partial([uv])=e_v-e_u$. 
That is, $\partial(z)= \sum_{[uv]\in E}z_{[uv]} (e_v-e_u)$ for a chain $z=\sum_{[uv]\in E}z_{[uv]} [uv]$.
A chain $z$ is a \emph{cycle}
if $\partial(z)=0$, or equivalently, if at every vertex the sum of
the coefficients on incoming edges equals the sum on outgoing edges.
Here a cycle need not be a single simple graph cycle.
It may be any real linear combination of oriented simple graph cycles.

\begin{example}
Let $\vec{G}$ be an oriented triangle on vertices $a,b,c$ with edges $(a,b)$, $(c,a)$, and $(b,c)$.
Then, $C_1(\vec{G};\mathbb{R})$ is isomorphic to $\mathbb{R}^3$, and the chains are of the form $z = \alpha[ab] + \beta[ca] + \gamma[bc]$ for $\alpha,\beta,\gamma \in \mathbb{R}$.
Setting $\alpha=\beta=\gamma$ gives a cycle, since we then have $\partial(z) = \alpha\partial([ab]) + \beta\partial([ca]) + \gamma\partial([bc]) = \alpha(e_b-e_a) + \beta(e_a - e_c) + \gamma(e_c - e_b) = \alpha(e_b-e_a + e_a - e_c + e_c - e_b) = 0$.
\end{example}

For distinct edges $e$ and $f$ of $\vec{G}$ that share a crossing $c$ in $D$, we use $c_{e,f}$ to denote the pair $(c,(e,f))$.
The \emph{sign of} $c_{e,f}$ is the sign of the determinant $\det(t_e,t_f)$ where $t_e$ and $t_f$ are vectors determined by the line segments forming $c_{e,f}$ pointing out in the direction of $e$ and $f$, respectively. 
Geometrically, the sign of $c_{e,f}$ is positive if and only if $f$ goes from the right side of $e$ to the left side of $e$.
Note that the sign of $c_{e,f}$ is opposite to the sign of $c_{f,e}$.

For distinct edges $e$ and $f$ of $\vec{G}$, let $I_D(e,f)$ be the sum of the signs of pairs $c_{e,f}$ taken over all crossings $c$ of $e$ and $f$ in $D$.
We set $I_D(e,e) = 0$ for every edge $e$ of $\vec{G}$ and note that $I_D(e,f)$ is the sign of $c_{e,f}$ if $D$ is simple and $e$ and $f$ cross at $c$.
In particular, $I_D(e,f) = 0$ if $e$ and $f$ do not cross in $D$.

For every vertex $v$ of $\vec{G}$, we fix a linear order $\prec_v$ of the edges incident to $v$ so that when these edges are traced counterclockwise in a small neighborhood of $v$, then we get the ordering $\prec_v$ after choosing an arbitrary gap in this radial counterclockwise order.
For edges $e$ and $f$ incident to a vertex $v$, we set 
\[
T_D^v(e,f)=
\begin{cases}
  1,  & e \prec_v f,\\
  -1, & f \prec_v e,\\
  0,  & e=f,
\end{cases}
\;\;\; \text{ and } \;\;\;
\sigma_D(v,e)=
\begin{cases}
  1,  & e \text{ points out of } v,\\
  -1, & e \text{ points into } v.
\end{cases}
\]
Observe that $\sigma_D(v, \bar e) = -\sigma_D(v,e)$ for the reverse orientation $\bar e$ of $e$.

Finally, for two not necessarily distinct edges $e$ and $f$ of $\vec{G}$, we define
\begin{equation}
\label{eq-J}
J_D(e,f)=I_D(e,f)+\frac{1}{2}
\sum_{v \in V : v \in e \cap f}
\sigma_D(v,e)\sigma_D(v,f)T_D^v(e,f).
\end{equation}
Note that $J_D(e,e) = 0$ for every edge $e$ of $\vec{G}$ and that $J_D(e,f) = -J_D(f,e)$ for all edges $e$ and $f$ of $\vec{G}$.
The choice of $\sigma_D(v,e)$ guarantees that $J_D(\bar e,f) = - J_D(e,f)$ and $J_D(e,\bar f) = -J_D(e,f)$ for reverse orientations $\bar e$ and $\bar f$ of $e$ and $f$, respectively. 
Moreover, the second term in $J_D(e,f)$ is chosen so that if we replace a crossing $c$ at $(e,f)$ with a vertex $v$, then the sign of $c_{e,f}$ corresponds to the value that $v$ adds in the second term over its four incident alternating spokes determined by portions of $e$ and $f$.
We extend $J_D$ to arbitrary pairs $(z,z')$ of chains $z=\sum_e z_e e$ and $z' = \sum_f z'_f f$ from $C_1(\vec{G};\mathbb{R})$ by letting $J_D(z,z') = \sum_{e,f} z_ez'_f J_D(e,f)$, obtaining a bilinear extension of~\eqref{eq-J}.
We extend $I_D$ bilinearly in the same way.

We now prove a crucial property of the bilinear extension of $J_D$.

\begin{lemma}
\label{lem-cycles}
Let $\vec{G}$ be an oriented graph and $D$ be a polygonal drawing of $\vec{G}$.
Any two cycles $z$ and $z'$ from $C_1(\vec{G};\mathbb{R})$ satisfy $J_D(z,z')=0$. 
\end{lemma}
\begin{proof}
Let $\widehat D$ be the plane drawing of a multigraph $\widehat{G}$
obtained from $D$ by replacing every crossing by a vertex.
For each edge $e$ of $\vec{G}$, we orient each subdivided piece of $e$ consistently with $e$.
If $z=\sum_{e \in E} z_e e$, then we assign each piece of $e$ the coefficient $z_e$, lifting $z$ to $\widehat z \in C_1(\widehat{G};\mathbb{R})$.
We proceed analogously for $z'$, lifting it to $\widehat z' \in C_1(\widehat{G};\mathbb{R})$.
The boundary equations remain unchanged at the original vertices, and the incoming and outgoing pieces of each subdivided edge cancel at every new vertex.
Thus, $\widehat z$ and $\widehat z'$ are cycles of $\widehat{G}$.

We apply the definition of $J_D$ also to $\widehat D$, retaining the chosen cuts at radial orders around the original vertices and choosing arbitrary cuts at the new ones. 
Note that the choice of cut at a new vertex $v$ does not change the contribution of $v$ to $J_{\widehat D}(\widehat z,\widehat z')$, since moving the edge $h$ that is the first incident to $v$ in $\prec_v$ to the last position changes the contribution by
\[
\sigma_{\widehat D}(v,h)\left(
\widehat z'_h\sum_{f\colon v\in f}\sigma_{\widehat D}(v,f)\widehat z_f
-
\widehat z_h\sum_{f\colon v\in f}\sigma_{\widehat D}(v,f)\widehat z'_f
\right)=0,
\]
where both sums vanish because $\widehat z$ and $\widehat z'$ are cycles.

We now show that $
J_D(z,z')=J_{\widehat D}(\widehat z,\widehat z')$.
Since $\widehat D$ has no crossings, $I_{\widehat D}(\widehat z,\widehat z')$ is zero and thus $J_{\widehat D}(\widehat z,\widehat z')$ consists only of vertex contributions.
At the original vertices of $\vec{G}$, these agree with those in $J_D(z,z')$.
At a crossing $c$ of edges $e,f$, with the sign $\varepsilon_{c_{e,f}}$ of $c_{e,f}$,
the four alternating spokes give the contribution
$
\varepsilon_{c_{e,f}}\bigl(z_ez'_f-z_fz'_e\bigr)
$
to $J_{\widehat D}(\widehat z,\widehat z')$.
Indeed, when $\varepsilon_{c_{e,f}}=1$, we can choose the cut just before the outgoing spoke of $e$, as this does not change the contribution to $J_{\widehat D}(\widehat z,\widehat z')$.
The four spokes then occur in the order $e_{\mathrm{out}},f_{\mathrm{out}},e_{\mathrm{in}},f_{\mathrm{in}}$.
Their orientation signs $\sigma_{\widehat D}$ are $+1,+1,-1,-1$, so the four mixed terms in~\eqref{eq-J} give $\frac12(1-1+1+1)=1$ as the coefficient of $z_ez'_f$ and, by antisymmetry, $-1$ as the coefficient of $z_fz'_e$.
Both signs reverse when $\varepsilon_{c_{e,f}}=-1$.
The expression $\varepsilon_{c_{e,f}}\bigl(z_ez'_f-z_fz'_e\bigr)$ is exactly the total contribution of the crossing $c$ to $I_D(z,z')$.
Consequently,
$
J_D(z,z')=J_{\widehat D}(\widehat z,\widehat z')$ and thus it suffices to show $J_{\widehat D}(\widehat z,\widehat z') = 0$.

For each bounded face $F$ of $\widehat D$, let $c_F$ be
the sum of its oriented boundary walks, counted with
multiplicity and oriented with $F$ on their left.
Let $\vec{C}$ be a simple cycle in $\widehat{G}$ oriented counterclockwise (possibly of length 2). 
As an element of $C_1(\widehat{G};\mathbb{R})$, we have $\vec{C} = \sum_{e \in E(\widehat{G})} c_e e$, where $c_e=1$ if the orientation of $e$ agrees with the orientation of $\vec{C}$, $c_e=-1$ if it does not, and $c_e=0$ if $e$ is not an edge of $\vec{C}$.
Then, 
$
\vec{C}=\sum_{F\text{ inside }\vec{C}}c_F.
$
Indeed, the contributions of each interior edge cancel, since
the facial boundaries traverse it in opposite directions,
whereas each edge of $\vec{C}$ occurs once, with the orientation of $\vec{C}$.
Reversing the orientation of $\vec{C}$ changes the sign of this sum.
Since every cycle is a real linear combination of simple cycles,
the oriented facial boundary chains span the cycle space ${\rm Ker}\, \partial$.
We may therefore write
$
\widehat z=\sum_F\lambda_F c_F$, where
$\lambda_F\in\mathbb R$ and the sum ranges over the bounded faces.
By bilinearity, $J_{\widehat D}(\widehat z,\widehat z')=\sum_F\lambda_F J_{\widehat D}(c_F,\widehat z')$.
Thus, it suffices to prove that
$J_{\widehat D}(c_F,\widehat z')=0$ for every bounded face $F$.
We may therefore assume that $\widehat z$ is the facial
boundary $c_F$ of $F$, oriented with the face $F$ on its left.

Writing $\widehat z' = \sum_{e \in E(\widehat{G})} \widehat z'_e e$, we set $b_v(e)=\sigma_{\widehat D}(v,e)\widehat z'_e$ for a vertex $v$ incident to $e$.
That is, $b_v(e) = \widehat z'_e$ if $e$ goes out of $v$ and $b_v(e) = -\widehat z'_e$ otherwise.
Since $\widehat z'$ is a cycle, we have
$\sum_{e \colon v \in e}b_v(e)=0$ at every vertex $v$.
Consider a corner $\kappa$ of $F$ at $v$ bounded by consecutive incident edges $e$ and $g$ in the counterclockwise order around $v$.
For each edge $f$ incident to $v$, we have
\[
T_{\widehat D}^v(e,f)-T_{\widehat D}^v(g,f)
=
\mathbf1_{f=e}+\mathbf1_{f=g}+c,
\qquad
c=
\begin{cases}
0, & \text{if the cut lies outside }\kappa,\\
-2, & \text{if the cut lies in }\kappa,
\end{cases}
\]
where $\mathbf1_A$ is $1$ if $A$ holds and $0$ otherwise.
Indeed, when $e\ne g$ and the cut lies outside $\kappa$,
we have $e\prec_v g$, and the difference is $1$ for
$f\in\{e,g\}$ and $0$ otherwise.
When $e\ne g$ and the cut lies in $\kappa$, we have $g\prec_v e$,
and the difference is $-1$ for $f\in\{e,g\}$ and $-2$ otherwise.
If $e=g$, then $v$ has degree one, the cut lies in $\kappa$,
and $f=e=g$, so the identity reduces to $0=1+1-2$.

By the definition of $J_{\widehat D}(c_F,\widehat z')$, the vertex $v$ contributes
\[
\frac{1}{2} \sum_{h \colon v \in h}\sum_{f \colon v \in f}
\sigma_{\widehat D}(v,h) \widehat z_h
\sigma_{\widehat D}(v,f) \widehat{z'}_f
T^v_{\widehat D}(h,f)
\]
to $J_{\widehat D}(c_F,\widehat z')$.
Since $\widehat z=c_F$, the coefficients
$\sigma_{\widehat D}(v,h)\widehat z_h$ are obtained by summing
the contributions of all corners of $F$ at $v$.
The boundary traversal with $F$ on its left enters $v$ along $g$
and leaves $v$ along $e$.
Hence, the corner $\kappa$ contributes
$\mathbf1_{h=e}-\mathbf1_{h=g}$
to $\sigma_{\widehat D}(v,h)\widehat z_h$
for each incident edge $h$, including zero when $e=g$.
Thus, the corner $\kappa$ contributes
\[
\frac{1}{2}\sum_{f\colon v\in f}
\bigl(T_{\widehat D}^v(e,f)-T_{\widehat D}^v(g,f)\bigr)b_v(f)
=\frac12\bigl(b_v(e)+b_v(g)\bigr)
\]
to $J_{\widehat D}(c_F,\widehat z')$, where in the case $c=-2$ we use $\sum_{f\colon v\in f}b_v(f)=0$, which holds as $\widehat z'$ is a cycle.
Summing over all corners, each occurrence of a boundary edge
$e=(v,w)$ contributes $\frac12b_v(e)+\frac12b_w(e)=0$, since $\sigma_{\widehat D}(v,e)=-\sigma_{\widehat D}(w,e)$.
This cancellation also applies when the facial boundary repeats
vertices or edges. Hence
$J_{\widehat D}(\widehat z,\widehat z')=0$ for every facial boundary,
and therefore for every cycle $\widehat z$.
It follows that $J_D(z,z')=0$.
\end{proof}

\subsection{Lower bound on the number of crossings}

In this subsection, we finish the proof of Theorem~\ref{thm-HararyHill} by estimating the number of crossings in every drawing of $K_n$ from below by $Z(n)$.
Let $D$ be a polygonal drawing of an oriented $\vec{K}_n$.
We associate each vertex $v_i$ of $\vec{K}_n$ with a real number $t_i$ so that the numbers $t_1,\dots,t_n$ are pairwise distinct.
We first use the fact that the bilinear function $J_D$ is zero on any pair of cycles to show that we can choose suitable weights on pairs of edges of $K_n$ so that, for every polynomial $H(u,v,x,y)$ of degree at most $n-2$ in each variable, the resulting weighted sum over all 4-tuples of vertices vanishes.

To define the weights, we set $F(x) = \prod_{i=1}^n (x-t_i)$ and $w_i = 1/F'(t_i)$, where $F'$ is the derivative of $F$.
Note that $F'(t_i) = \prod_{j \neq i} (t_i-t_j)$, so $w_i = 1/F'(t_i) = 1/\prod_{j \neq i} (t_i-t_j)$ is properly defined.
The weights $w_i$ are chosen so that $\sum_{i=1}^n w_it_i^a = 0$ for every integer $a$ with $0 \leq a \leq n-2$.
For a positive integer $n$, we write $[n]$ to denote $\{1,\dots,n\}$.

\begin{lemma}
\label{lem-4.1}
Let $n \geq 2$ be an integer and let $D$ be a polygonal drawing of an oriented $\vec{K}_n$.
Then, every polynomial $H(u,v,x,y)$ of degree at most $n-2$ separately in each variable satisfies 
\begin{equation}
\label{eq-lem4.1}
\sum_{i,j,k,l \in [n]} w_iw_jw_kw_l\cdot  J_D([v_iv_j],[v_kv_l]) \cdot H(t_i,t_j,t_k,t_l) =0.
\end{equation}
\end{lemma}
\begin{proof}
For integers $0\leq a,b\leq n-2$, we define the chain $z_{a,b}=\sum_{i,j=1}^nw_iw_j t_i^a t_j^b[v_iv_j]$.
Note that this is indeed a chain as $\vec{K}_n$ contains all edges.
The boundary of $z_{a,b}$ is
\[
\partial(z_{a,b})
= \sum_{i,j=1}^n w_iw_j t_i^a t_j^b
   (e_{v_j}-e_{v_i})
= \sum_{k=1}^n w_k
   \left(
   t_k^b\sum_{i=1}^n w_i t_i^a
   -t_k^a\sum_{j=1}^n w_j t_j^b
   \right)e_{v_k}.
\]
The coefficient $w_k\left(t_k^b\sum_i w_i t_i^a-t_k^a\sum_j w_j t_j^b\right)$ at $e_{v_k}$ is zero by the choice of the weights, as Lemma~\ref{lem-3.1} gives $\sum_{i=1}^n w_it_i^c = 0$ for every integer $c$ with $0 \leq c \leq n-2$.
Thus, $\partial(z_{a,b})=0$ and every $z_{a,b}$ is a cycle.

If $H$ is a monomial of the form $H(u,v,x,y)=u^a v^b x^c y^d$ with
$0\leq a,b,c,d\leq n-2$, then, using bilinearity of $J_D$, the left-hand side of~\eqref{eq-lem4.1} becomes $J_D(z_{a,b},z_{c,d})$, which is zero by Lemma~\ref{lem-cycles}, as $z_{a,b}$ and $z_{c,d}$ are cycles.
In general, the polynomial $H$ of degree at most $n-2$ separately in each
variable is a linear combination of these monomials, so~\eqref{eq-lem4.1} holds by linearity.
\end{proof}

We can now proceed with the proof of the Harary--Hill conjecture.

\begin{proof}[Proof of Theorem~\ref{thm-HararyHill}]
Since there are known constructions of drawings of $K_n$ with $Z(n)$ crossings for every positive integer $n$~\cite{BlazekKoman1964,DamianiDAntonaSalemi1994,HararyHill1963}, it suffices to prove the tight lower bound.
Let $D$ be a drawing of $K_n$ with a minimum number of crossings.
Without loss of generality, we may assume that $D$ is simple and polygonal~\cite{Guy1972}.
First, the statement is trivial for $n \leq 4$ as $Z(n) = 0$, so we assume that $n \geq 5$.
It also suffices to prove the conjecture for odd values of $n$, as the cases with even $n$ then follow by a standard double-counting argument.
Indeed, assuming ${\rm cr}(K_{2s+1}) \geq Z(2s+1) = \binom{s}{2}^2$ for $s\geq 2$, double-counting the pairs $(c,v)$ where $v$ is a vertex not incident to a crossing $c$ in a simple drawing $D$ of $K_{2s+2}$, with $s\geq 2$, gives
\[(2s-2){\rm cr}(D) \geq (2s+2)\binom{s}{2}^2,\]
so ${\rm cr}(D) \geq (s+1)s^2(s-1)/4=Z(2s+2)$.
Thus, in the rest of the proof, we assume that $n=2s+1$ for an integer $s\ge2$.

Let the vertices of $K_n$ be $v_1,\dots,v_n$ and fix an arbitrary orientation $\vec{K}_n$ of $K_n$.
We associate each vertex $v_i$ with a real number $t_i$ so that the numbers $t_1,\dots,t_n$ are distinct.
Consider the bilinear extension $J_D \colon C_1(\vec{K}_n;\mathbb{R}) \times C_1(\vec{K}_n;\mathbb{R}) \to \mathbb{R}$ of the function from~\eqref{eq-J}.
We again set $F(t) = \prod_{i=1}^n (t-t_i)$ and $w_i = 1/F'(t_i)$.

Let $\mathcal{X}_{\pm}$ be the set of ordered pairs of undirected edges that cross in $D$.
Since $D$ is simple, no pair in $\mathcal{X}_{\pm}$ shares an endpoint.
Among pairs of edges with no common endpoint, these are exactly the pairs $(e,f)$ with $J_D(e,f)\neq0$ in any orientation of $K_n$.
By definition, we have ${\rm cr}(D) = |\mathcal{X}_{\pm}|/2$.

Let $\mathcal{S}$ be the vector space of polynomials $P(u,v,x,y)$ of degree at most $s-2$ separately in each variable that are symmetric within $(u,v)$ and $(x,y)$.
That is, $P(u,v,x,y) = P(v,u,x,y) = P(u,v,y,x) = P(v,u,y,x)$.
The polynomials $u^av^b+u^bv^a$, with $0\le a\le b\le s-2$,
form a basis of the space of symmetric bivariate polynomials of degree at most $s-2$ in each variable.
There are $\binom{s-1}{2} + s-1 = \binom{s}{2}$ options for such polynomials.
Thus, taking products of these basis elements in $(u,v)$ and $(x,y)$
gives a basis of $\mathcal S$, so $\dim\mathcal S = \binom{s}{2}^2$.

For $P\in\mathcal S$ and $\alpha=(\{v_i,v_j\},\{v_k,v_l\})\in\mathcal{X}_{\pm}$,
write $P_\alpha=P(t_i,t_j,t_k,t_l)$.
By the symmetry of $P$, this value does not depend on how we list
the endpoints of either edge.
We define the polynomial
\[
\Delta(u,v,x,y)=(u-v)(x-y)(u-x)(u-y)(v-x)(v-y),
\]
which is used to eliminate contributions from edges sharing an endpoint and to compensate for the sign changes of $J_D$ under edge reversal.
If the two edges $(u,v)$ and $(x,y)$ share an endpoint, then one of the last four factors of $\Delta$ is zero, so $\Delta=0$.
If we reverse either edge, then the corresponding factor among $(u-v)$ and $(x-y)$ changes sign, while the product of the last four factors remains unchanged.
Thus, $\Delta$ changes sign as well. 
Note that $\Delta$ is of degree at most 3 separately in each coordinate.

For any $P\in\mathcal S$, the polynomial $H = P^2\Delta$ has degree at most $2(s-2)+3=2s-1$ separately in each variable.
Thus, by Lemma~\ref{lem-4.1}, 
\begin{equation}
\label{eq-HJ}
\sum_{i,j,k,l \in [n]} w_iw_jw_kw_l\cdot  J_D([v_iv_j],[v_kv_l]) \cdot H(t_i,t_j,t_k,t_l) =0,
\end{equation}
where nonzero terms in the sum come only from $\mathcal{X}_{\pm}$, as the factor $\Delta$ eliminates repeated endpoints, and $J_D$ vanishes on pairs of noncrossing edges with no common endpoint.

For each ordered pair in $\mathcal{X}_{\pm}$, its four choices of orientations give equal contributions to~\eqref{eq-HJ}.
Reversing either edge changes the signs of both $J_D$ and $\Delta$, while leaving $P^2$ and the product $w_iw_jw_kw_l$ unchanged.
Thus, after division by four, \eqref{eq-HJ} gives $\sum_{\alpha\in\mathcal{X}_{\pm}} a_\alpha P_\alpha^2=0$ where
\[
a_{(\{v_i,v_j\},\{v_k,v_l\})}
=w_iw_jw_kw_l\cdot J_D([v_iv_j],[v_kv_l]) \cdot \Delta(t_i,t_j,t_k,t_l),
\]
and every $a_\alpha$ is nonzero.
For $\alpha=(\{v_i,v_j\},\{v_k,v_l\})$, exchanging the two edges gives a different ordered pair $\alpha'=(\{v_k,v_l\},\{v_i,v_j\})$ with $a_\alpha=-a_{\alpha'}$, since this exchange preserves $\Delta$ and negates $J_D$.
Hence the set $\mathcal{X}=\{\alpha\in\mathcal{X}_{\pm}:a_\alpha>0\}$ has size $|\mathcal{X}_{\pm}|/2={\rm cr}(D)$.

Suppose for contradiction that ${\rm cr}(D)<\dim\mathcal S$.
Since $|\mathcal{X}|={\rm cr}(D)$, the conditions $P_\alpha=0$ for $\alpha\in\mathcal{X}$ with variables corresponding to the coefficients of $P$ form fewer
homogeneous linear equations than $\dim\mathcal S$, so they have a nonzero solution $P\in\mathcal S$.
Choose such a $P$ of minimum degree in $u$.
The terms with $a_\alpha>0$ in $\sum_{\alpha\in\mathcal{X}_{\pm}}a_\alpha P_\alpha^2=0$
vanish, and all remaining coefficients are negative.
Since $P_\alpha^2\ge0$, every remaining term must also vanish.
Thus, $P_\alpha=0$ for every $\alpha\in\mathcal{X}_{\pm}$.

We fix a vertex $v_r$ of $\vec{K}_n$ and define $F_r(x)=\prod_{j\neq r}(x-t_j)$. We factor the polynomial $F_r(x)$ as
$
F_r(x)=U(x)V(x),
$
where $\deg U=\deg V=s$.
This can be done as the number of linear factors in $F_r(t)$ is $n-1=2s$.
Recall that $w_r = 1/F'(t_r) = 1/F_r(t_r)$.

For any $Q\in\mathcal P_{s-1,s-1}$, we define the polynomial
\[
\bar H(u,v,x,y)=P(u,v,x,y)(v-y)(v-x)(y-u)U(u)V(x)Q(v,y).
\]
Note that its degree in each variable is at most $2s-1=n-2$.
Thus, Lemma~\ref{lem-4.1} applied to $\bar H$ gives 
\begin{equation}
\label{eq-claim1Eq1}
\sum_{i,j,k,l \in [n]} w_iw_jw_kw_l\cdot J_D([v_iv_j],[v_kv_l]) \cdot \bar H(t_i,t_j,t_k,t_l)=0.
\end{equation}
The terms with $i=j$ or $k=l$ vanish because $[v_iv_i]=[v_kv_k]=0$. 
The terms corresponding to a pair of disjoint edges vanish because $P_\alpha=0$ for $\alpha\in\mathcal{X}_{\pm}$ and $J_D=0$ for all other such pairs. 
The factor $(v-y)(v-x)(y-u)$ eliminates
all terms corresponding to edges with a shared vertex, except the case $u=x$. Moreover,
$U(t_j)V(t_j)=F_r(t_j)=0$ for $j\neq r$, so only $u=x=t_r$ remains.
So after substituting $(u,v,x,y) = (t_r,t_j,t_r,t_l)$ and expanding the definition of $\bar H$, \eqref{eq-claim1Eq1} becomes
\[
w_r^2U(t_r)V(t_r)
\sum_{\substack{j\ne r\\l\ne r}}
J_D([v_rv_j],[v_rv_l])\,w_jw_l\cdot
P(t_r,t_j,t_r,t_l)
(t_j-t_l)(t_j-t_r)(t_l-t_r)Q(t_j,t_l)=0.
\]

Since adjacent edges do not cross in $D$, we get $J_D([v_rv_j],[v_rv_l])=\frac{1}{2}T_D^{v_r}([v_rv_j],[v_rv_l])$.
Using $(t_j-t_r)w_j=\frac{1}{F_r'(t_j)}$ for $j\neq r$, $w_r^2U(t_r)V(t_r)=w_r\neq0$,
and dividing the above equation by $w_r$, we obtain
\[
\sum_{\substack{j\neq r\\l\neq r}}
\frac{T_D^{v_r}([v_rv_j],[v_rv_l])}{2}\frac{t_j-t_l}{F_r'(t_j)F_r'(t_l)}
P(t_r,t_j,t_r,t_l)Q(t_j,t_l)=0.
\]
The left-hand side is the expression from Lemma~\ref{lem-3.3} on the $2s$ nodes $t_j$ with $j\neq r$, with $t_j\prec t_l$ if $[v_rv_j]\prec_{v_r}[v_rv_l]$.
Since $Q$ was arbitrary and $P(t_r,v,t_r,y)$, viewed as a polynomial in $v,y$, belongs to $\mathcal P_{s-2,s-2}$, Lemma~\ref{lem-3.3} gives
$P(t_r,v,t_r,y)\equiv0$ for every $r$.

The polynomial $P(u,v,u,y)$ has degree at most $2s-4<n$ in $u$ and vanishes at all $n$ distinct labels $t_r$.
Each of its coefficients as a polynomial in $v,y$ therefore vanishes identically in $u$.
Thus, $P(u,v,u,y)\equiv0$, so $u-x$ divides $P$.
The two within-pair symmetries imply that $u-y$, $v-x$, and $v-y$ also divide $P$.
Since these four linear factors are pairwise non-proportional, their product $C=(u-x)(u-y)(v-x)(v-y)$ divides $P$.
Since $C$ has the same two symmetries and is nonzero at every
evaluated pair of disjoint edges, $P/C$ is a nonzero element
of $\mathcal S$ with $(P/C)_\alpha=0$ for every $\alpha\in\mathcal{X}$. Its degree in $u$ is two
less than that of $P$, contradicting the choice of $P$.
This contradiction gives 
\[
{\rm cr}(D) \ge\dim\mathcal S=\binom{s}{2}^{\!2}=Z(2s+1),
\]
which finishes the proof.
\end{proof}

\paragraph{Comparison with the original proof.}
Our proof simplifies several steps of the original argument~\cite{openAI}.
We establish the cycle identity in Lemma~\ref{lem-cycles} using planarization and facial boundaries, avoiding any perturbation argument.
We also prove a weaker version of Lemma~\ref{lem-3.3} sufficient for our purposes, and isolate individual vertices directly by factoring $F_r=UV$, avoiding an additional interpolation argument.
Most importantly, applying the cycle identity to $P^2\Delta$ yields a direct proof by contradiction based on the signs of the coefficients, replacing the original argument involving a nondegenerate bilinear form and orthogonal complements.

\subsection{Other crossing number variants}
\label{subsec-variants}

The proof of Theorem~\ref{thm-HararyHill} also establishes
${\rm ocr}_+(K_n)={\rm pcr}_+(K_n)=Z(n)$.
Indeed, let $D$ be an arbitrary semisimple drawing of $K_n$, and consider two independent edges
$e=ab$ and $f=cd$.
Since adjacent edges do not cross, the triangle $abc$
forms a simple closed curve, and $f$ meets its boundary
only at $c$ and at its crossings with~$e$.
As $f$ is traversed from $c$ to $d$, successive crossings
with $e$ alternate between entering and leaving the region
bounded by the triangle.
Consequently, their signs alternate, and hence
$I_D(e,f)\in\{-1,0,1\}$, with $I_D(e,f)\neq0$ if and only if
$c_D(e,f)$ is odd.

We can therefore repeat the proof of Theorem~\ref{thm-HararyHill}
without assuming that $D$ is simple, defining
$\mathcal{X}_{\pm}$ as the set of ordered pairs of independent
edges $(e,f)$ with $I_D(e,f)\neq0$.
Indeed, for such edges $J_D(e,f)=I_D(e,f)$, so all the
algebraic identities used in the proof remain valid.
Moreover, $|\mathcal{X}_{\pm}|/2$ is precisely the number
of pairs of edges crossing oddly in $D$.
The same argument thus gives
${\rm ocr}_+(K_n)\ge Z(n)$ for odd $n$.
For even $n$, the same double-counting argument applies, since each pair of oddly crossing edges survives exactly
$n-4$ vertex deletions.
Since
\[
{\rm ocr}_+(K_n)\le {\rm pcr}_+(K_n)
\le {\rm cr}(K_n)=Z(n),
\]
where the last equality follows from Theorem~\ref{thm-HararyHill},
both variants attain the claimed value.

\section{Immediate applications}
\label{sec-applications}

We combine Theorem~\ref{thm-HararyHill} with known reductions
to obtain exact formulas for crossing numbers of several other families of graphs that were considered in the literature.

\subsection{Complete graphs with deleted edges}

The following consequence of Theorem~\ref{thm-HararyHill} follows from
Felsner's reduction~\cite[Theorem~9]{Felsner2026}.
We include the short argument for the lower bound.

\begin{cor}
\label{cor-perfectMatching}
For every integer $k\geq2$ and every perfect matching $M$ of $K_{2k}$,
\[
{\rm cr}(K_{2k}-M)=6\binom{k}{4}.
\]
\end{cor}

\begin{proof}
The upper bound is given by the construction discussed
in~\cite[Section~3]{Felsner2026}.
For the lower bound, let $D$ be an optimal simple polygonal drawing
of $K_{2k}-M$, and choose a vertex~$v$.
For each edge $xy\in M$ not incident to~$v$, join $x$ to $y$ along
the path $xvy$, and then delete $v$.
Outside a small disk around $v$, each segment of the original drawing
is used at most once, so the number of crossings does not increase there.
Inside the disk, join the corresponding boundary points by chords
and perturb them so that no three meet.
These $k-1$ chords create at most $\binom{k-1}{2}$ crossings.
The resulting drawing $D'$ after removing $v$ is a drawing of $K_{2k-1}$.
Thus, Theorem~\ref{thm-HararyHill} gives
\[
{\rm cr}(D)
\geq {\rm cr}(D')-\binom{k-1}{2}
\geq Z(2k-1)-\binom{k-1}{2}
=6\binom{k}{4},
\]
finishing the proof.
\end{proof}

For a single removed edge, we use a reduction of Ouyang, Wang, and Huang~\cite{OuyangWangHuang2014}.

\begin{cor}
\label{cor-deletedEdge}
For every integer $k\geq2$ and every edge $e$ of $K_{2k}$,
\[
{\rm cr}(K_{2k}-e)
=Z(2k)-\binom{k-1}{2}
=\frac{(k^2-1)(k-2)^2}{4}.
\]
\end{cor}

\begin{proof}
Ouyang, Wang, and Huang~\cite[Theorem~3.2]{OuyangWangHuang2014} showed that ${\rm cr}(K_{2k}-e)
\geq\frac{2k+2}{2k-2}{\rm cr}(K_{2k-1})$ for every integer $k \geq 2$. 
Together with Theorem~\ref{thm-HararyHill}, this gives
\[
{\rm cr}(K_{2k}-e)
\geq\frac{2k+2}{2k-2}{\rm cr}(K_{2k-1})
=\frac{k+1}{k-1}Z(2k-1)
=\frac{(k^2-1)(k-2)^2}{4}.
\]
The matching upper bound is given
in~\cite[Lemma~2.3]{OuyangWangHuang2014}.
\end{proof}

\subsection{Cartesian products of complete graphs with paths}

Let $P_t$ denote the path on $t$ vertices.
The \emph{Cartesian product} $G \square H$ of two graphs $G$ and $H$
has vertex set $V(G)\times V(H)$, where two vertices $(u,v)$ and $(u',v')$ are adjacent if and only if
$u=u'$ and $vv'\in E(H)$, or
$v=v'$ and $uu'\in E(G)$.
The next consequence uses a result of Ouyang, Wang, and
Huang~\cite{OuyangWangHuang2014} and determines the crossing number of $K_{2k}\square P_t$.

\begin{cor}
\label{cor-cliquePath}
For all integers $k\geq2$ and $t\geq2$,
\[
{\rm cr}(K_{2k}\square P_t)
=\frac{k(k-1)^2}{4}\bigl(t(k+2)-4\bigr).
\]
\end{cor}

\begin{proof}
In~\cite{OuyangWangHuang2014}, $P_n$ denotes a path with $n+1$ vertices, so we use their formula with $n=t-1$.
The decomposition proved by Ouyang, Wang, and
Huang~\cite[Corollary~5.1]{OuyangWangHuang2014} yields
\[
{\rm cr}(K_{2k}\square P_t)
=(t-2){\rm cr}(K_{2k+2}-e)+2{\rm cr}(K_{2k+1}),
\]
where $e$ is any edge of $K_{2k+2}$.
By Corollary~\ref{cor-deletedEdge} and
Theorem~\ref{thm-HararyHill}, the right-hand side equals
\[
(t-2)\frac{k(k+2)(k-1)^2}{4}
+2\frac{k^2(k-1)^2}{4}
=\frac{k(k-1)^2}{4}\bigl(t(k+2)-4\bigr).
\]
\end{proof}

\subsection{Cartesian products of complete graphs with even cycles}

Using a result of Zheng, Lin, Yang, and Deng~\cite{ZhengLinYangDeng2008}, we can also determine the crossing numbers of Cartesian products of complete graphs with even cycles.

\begin{cor}
\label{cor-evenCycleProduct}
For all integers $m\geq5$ and all even integers $n\geq4$,
\[
{\rm cr}(K_m\square C_n)=nZ(m+2).
\]
\end{cor}

\begin{proof}
Zheng, Lin, Yang, and Deng~\cite[Theorems~2.2 and~3.1]{ZhengLinYangDeng2008}
proved that
\[
n\,{\rm cr}(K_{m+2})
\leq {\rm cr}(K_m\square C_n)
\leq nZ(m+2)
\]
for these parameters. 
Theorem~\ref{thm-HararyHill} then gives the equality.
\end{proof}

\subsection{Extended Sierpi\'nski graphs}

We can also determine the crossing number of the two
regularizations of the Sierpi\'nski graph $S(n,k)$ defined by
Klav\v{z}ar and Mohar~\cite[Section~2]{KlavzarMohar2005}.
For integers $n\geq1$ and $k\geq2$, the
\emph{Sierpi\'nski graph} $S(n,k)$ has vertex set $[k]^n$.
Two distinct vertices $u=(u_1,\dots,u_n)$ and
$v=(v_1,\dots,v_n)$ are adjacent if there exists
$h\in[n]$ such that
$
u_i=v_i$ for $i<h$, $u_h\neq v_h$, and $u_i=v_h,\ v_i=u_h$ for $i>h$.
The $k$ vertices $(i,\dots,i)$, $i\in[k]$, are called \emph{extreme vertices}.
The graph $S^+(n,k)$ is obtained from $S(n,k)$
by adding a new vertex adjacent to all its
extreme vertices.
For $n=1$, we set $S^{++}(1,k)=K_{k+1}$.
For $n\geq2$, the graph $S^{++}(n,k)$ is obtained
from the disjoint union of $S(n,k)$ and
$S(n-1,k)$ by adding a perfect matching between
their respective sets of extreme vertices.
Note that the graphs $S^{+}(n,k)$ and $S^{++}(n,k)$ are $k$-regular.

\begin{cor}
\label{cor-sierpinski}
For all integers $k\geq2$ and $n\geq1$,
\[
{\rm cr}(S^+(n,k))=\frac{k^n-1}{k-1}Z(k+1)
\]
and
\[
{\rm cr}(S^{++}(n,k))
=\frac{(k+1)k^{n-1}-2}{k-1}Z(k+1).
\]
\end{cor}

\begin{proof}
By~\cite[Theorem~4.1]{KlavzarMohar2005}, the two identities hold
with $Z(k+1)$ replaced by ${\rm cr}(K_{k+1})$.
The result therefore follows from Theorem~\ref{thm-HararyHill}.
\end{proof}

\paragraph{AI statement.}
ChatGPT (OpenAI, Astra model) was used to assist in understanding the original proof, exploring possible simplifications of some of its arguments, improving the clarity and presentation of the text, and identifying relevant literature. 
All mathematical arguments were independently checked, rewritten, and verified by the author.

\printbibliography[heading=bibintoc]

@article{Tutte1970,
  author  = {Tutte, W.T.},
  title   = {Toward a theory of crossing numbers},
  journal = {Journal of Combinatorial Theory},
  volume  = {8},
  number  = {1},
  pages   = {45--53},
  year    = {1970},
  doi     = {10.1016/S0021-9800(70)80007-2}
}

@article{Skopenkov2020,
  author  = {Skopenkov, A.},
  title   = {Invariants of graph drawings in the plane},
  journal = {Arnold Mathematical Journal},
  volume  = {6},
  pages   = {21--55},
  year    = {2020},
  doi     = {10.1007/s40598-019-00128-5}
}

@article{GoncharovKenyon2013,
  author  = {Goncharov, A.B. and Kenyon, R.},
  title   = {Dimers and cluster integrable systems},
  journal = {Annales scientifiques de l'\'{E}cole Normale
             Sup\'{e}rieure (4)},
  volume  = {46},
  number  = {5},
  pages   = {747--813},
  year    = {2013},
  doi     = {10.24033/asens.2201}
}

@misc{openAI,
  author = {{OpenAI}},
  title = {{The crossing number of complete graphs}},
  howpublished = {OpenAI Math Release preprint
                  \href{https://github.com/openai/math/blob/main/preprints/The-crossing-number-of-complete-graphs-September-23-2026/paper.pdf}{OAI:The-crossing-number-of-complete-graphs-September-23-2026}},
  year = {2026}
}

@article{HararyHill1963,
  author  = {Harary, Frank and Hill, Anthony},
  title   = {On the number of crossings in a complete graph},
  journal = {Proceedings of the Edinburgh Mathematical Society},
  volume  = {13},
  number  = {4},
  pages   = {333--338},
  year    = {1963},
  doi     = {10.1017/S0013091500025645}
}

@inproceedings{BlazekKoman1964,
  author    = {Bla{\v{z}}ek, J. and Koman, M.},
  title     = {A minimal problem concerning complete plane graphs},
  editor    = {Fiedler, M.},
  booktitle = {Theory of Graphs and Its Applications:
               Proceedings of the Symposium, Smolenice, 1963},
  pages     = {113--117},
  publisher = {Publishing House of the Czechoslovak Academy of Sciences},
  address   = {Prague},
  year      = {1964}
}

@article{DamianiDAntonaSalemi1994,
  author  = {Damiani, E. and D'Antona, O. and Salemi, P.},
  title   = {An upper bound to the crossing number of the complete
             graph drawn on the pages of a book},
  journal = {Journal of Combinatorics, Information \& System Sciences},
  volume  = {19},
  number  = {1--2},
  pages   = {75--84},
  year    = {1994}
}

@book{StoerBulirsch2002,
  author    = {Stoer, Josef and Bulirsch, Roland},
  title     = {Introduction to Numerical Analysis},
  edition   = {3},
  series    = {Texts in Applied Mathematics},
  volume    = {12},
  publisher = {Springer},
  address   = {New York},
  year      = {2002},
  doi       = {10.1007/978-0-387-21738-3}
}

@incollection{Guy1972,
  author    = {Guy, Richard K.},
  title     = {Crossing numbers of graphs},
  booktitle = {Graph Theory and Applications},
  editor    = {Alavi, Y. and Lick, D. R. and White, A. T.},
  series    = {Lecture Notes in Mathematics},
  volume    = {303},
  pages     = {111--124},
  publisher = {Springer},
  address   = {Berlin, Heidelberg},
  year      = {1972},
  doi       = {10.1007/BFb0067363}
}

@article{Dvir2009,
  author  = {Dvir, Zeev},
  title   = {On the size of {Kakeya} sets in finite fields},
  journal = {Journal of the American Mathematical Society},
  volume  = {22},
  number  = {4},
  pages   = {1093--1097},
  year    = {2009},
  doi     = {10.1090/S0894-0347-08-00607-3}
}

@misc{Felsner2026,
  author = {Felsner, Stefan},
  title = {Antipodal Pairs and Crossing Numbers of Complete Graphs},
  year = {2026},
  note = {Preprint},
  url = {https://page.math.tu-berlin.de/~felsner/Paper/antipodal.pdf},
  urldate = {2026-10-22}
}

@article{OuyangWangHuang2014,
  author = {Ouyang, ZhangDong and Wang, Jing and Huang, YuanQiu},
  title = {The crossing number of the {Cartesian} product of paths with complete graphs},
  journal = {Discrete Mathematics},
  volume = {328},
  pages = {71--78},
  year = {2014},
  doi = {10.1016/j.disc.2014.03.011}
}

@article{KlavzarMohar2005,
  author = {Klav{\v{z}}ar, Sandi and Mohar, Bojan},
  title = {Crossing numbers of {Sierpi{\'n}ski}-like graphs},
  journal = {Journal of Graph Theory},
  volume = {50},
  number = {3},
  pages = {186--198},
  year = {2005},
  doi = {10.1002/jgt.20107}
}

@article{ZhengLinYangDeng2008,
  author = {Zheng, Wenping and Lin, Xiaohui and Yang, Yuansheng and Deng, Chengrui},
  title = {On the crossing numbers of {$K_m\square C_n$} and {$K_{m,l}\square P_n$}},
  journal = {Discrete Applied Mathematics},
  volume = {156},
  number = {10},
  pages = {1892--1907},
  year = {2008},
  doi = {10.1016/j.dam.2007.09.007}
}

@article{Guy1960,
  author  = {Guy, Richard K.},
  title   = {A Combinatorial Problem},
  journal = {Nabla (Bulletin of the Malayan Mathematical Society)},
  volume  = {7},
  pages   = {68--72},
  year    = {1960}
}

@article{ErdosGuy1973,
  author  = {Erd{\H{o}}s, Paul and Guy, Richard K.},
  title   = {Crossing Number Problems},
  journal = {The American Mathematical Monthly},
  volume  = {80},
  number  = {1},
  pages   = {52--58},
  year    = {1973},
  doi     = {10.1080/00029890.1973.11993230}
}

@article{AbregoEtAl2013,
  author  = {{\'A}brego, Bernardo M. and
             Aichholzer, Oswin and
             Fern{\'a}ndez-Merchant, Silvia and
             Ramos, Pedro and Salazar, Gelasio},
  title   = {The 2-Page Crossing Number of {$K_n$}},
  journal = {Discrete \& Computational Geometry},
  volume  = {49},
  number  = {4},
  pages   = {747--777},
  year    = {2013},
  doi     = {10.1007/s00454-013-9514-0}
}

@article{AbregoEtAl2014,
  author  = {{\'A}brego, Bernardo M. and
             Aichholzer, Oswin and
             Fern{\'a}ndez-Merchant, Silvia and
             Ramos, Pedro and Salazar, Gelasio},
  title   = {Shellable Drawings and the Cylindrical
             Crossing Number of {$K_n$}},
  journal = {Discrete \& Computational Geometry},
  volume  = {52},
  number  = {4},
  pages   = {743--753},
  year    = {2014},
  doi     = {10.1007/s00454-014-9635-0}
}

@inproceedings{Aichholzer2021,
  author    = {Aichholzer, Oswin},
  title     = {Another Small but Long Step for Crossing Numbers:
               cr(13) = 225 and cr(14) = 315},
  booktitle = {Proceedings of the 33rd Canadian Conference
               on Computational Geometry (CCCG 2021)},
  pages     = {72--77},
  year      = {2021},
  url       = {https://cccg.ca/proceedings/2021/CCCG2021.pdf}
}

@article{BaloghLidickySalazar2019,
  author  = {Balogh, J{\'o}zsef and
             Lidick{\'y}, Bernard and Salazar, Gelasio},
  title   = {Closing in on Hill's Conjecture},
  journal = {SIAM Journal on Discrete Mathematics},
  volume  = {33},
  number  = {3},
  pages   = {1261--1276},
  year    = {2019},
  doi     = {10.1137/17M1158859}
}

@misc{BaloghEtAl2026,
  author = {Balogh, J{\'o}zsef and
            Layman, Nicholas and
            Lidick{\'y}, Bernard and
            Pfender, Florian and
            Salazar, Gelasio and
            Spiro, Sam and
            Volec, Jan},
  title  = {Crossing Numbers of Complete Bipartite Graphs},
  year   = {2026},
  eprint = {2610.25204},
  archivePrefix = {arXiv},
  primaryClass = {math.CO},
  url    = {https://arxiv.org/abs/2610.25204},
  note   = {Preprint}
}

@article{BalkoFulekKyncl2015,
  author  = {Balko, Martin and Fulek, Radoslav
             and Kyn{\v{c}}l, Jan},
  title   = {Crossing Numbers and Combinatorial
             Characterization of Monotone Drawings of {$K_n$}},
  journal = {Discrete \& Computational Geometry},
  volume  = {53},
  number  = {1},
  pages   = {107--143},
  year    = {2015},
  doi     = {10.1007/s00454-014-9644-z}
}

@article{AbregoEtAl2018,
  author  = {{\'A}brego, Bernardo M. and
             Aichholzer, Oswin and
             Fern{\'a}ndez-Merchant, Silvia and
             McQuillan, Dan and Mohar, Bojan and
             Mutzel, Petra and Ramos, Pedro and
             Richter, R. Bruce and Vogtenhuber, Birgit},
  title   = {Bishellable Drawings of {$K_n$}},
  journal = {SIAM Journal on Discrete Mathematics},
  volume  = {32},
  number  = {4},
  pages   = {2482--2492},
  year    = {2018},
  doi     = {10.1137/17M1147974}
}

@inproceedings{MutzelOettershagen2018,
  author    = {Mutzel, Petra and Oettershagen, Lutz},
  title     = {The Crossing Number of Seq-Shellable
               Drawings of Complete Graphs},
  booktitle = {Combinatorial Algorithms},
  series    = {Lecture Notes in Computer Science},
  volume    = {10979},
  pages     = {273--284},
  publisher = {Springer},
  year      = {2018},
  doi       = {10.1007/978-3-319-94667-2_23}
}

@inproceedings{StreltsovaWagner2025,
  author    = {Streltsova, Elizaveta and Wagner, Uli},
  title     = {Levels in Arrangements: Linear Relations,
               the {$g$}-Matrix, and Applications to
               Crossing Numbers},
  booktitle = {41st International Symposium on
               Computational Geometry (SoCG 2025)},
  series    = {Leibniz International Proceedings
               in Informatics (LIPIcs)},
  volume    = {332},
  pages     = {75:1--75:18},
  publisher = {Schloss Dagstuhl --
               Leibniz-Zentrum f{\"u}r Informatik},
  year      = {2025},
  doi       = {10.4230/LIPIcs.SoCG.2025.75}
}

@misc{gelasio,
  author = {Salazar, Gelasio},
  title = {The proof of the {H}arary--{H}ill conjecture: a detailed exposition and an extension to rotation systems},
  howpublished = {Arxiv preprint, arxiv:2611.03515},
  url          = {https://arxiv.org/abs/2611.03515},
  year = {2026}
}

@article{BeinekeWilson2010,
  author  = {Beineke, Lowell W. and Wilson, Robin J.},
  title   = {The Early History of the Brick Factory Problem},
  journal = {The Mathematical Intelligencer},
  volume  = {32},
  number  = {2},
  pages   = {41--48},
  year    = {2010},
  doi     = {10.1007/s00283-009-9120-4}
}

@article{FulekPelsmajerSchaeferStefankovic2012,
  author  = {Radoslav Fulek and Michael J. Pelsmajer
             and Marcus Schaefer and Daniel {\v S}tefankovi{\v c}},
  title   = {Adjacent Crossings Do Matter},
  journal = {Journal of Graph Algorithms and Applications},
  volume  = {16},
  number  = {3},
  pages   = {759--782},
  year    = {2012},
  doi     = {10.7155/jgaa.00266}
}

@misc{OpenAI2026Siegel,
  author       = {{OpenAI}},
  title        = {Uniform Exclusion of {Landau--Siegel} Zeros},
  year         = {2026},
  howpublished = {OpenAI Math Release preprint},
  note         = {October 1, 2026},
  url          = {https://github.com/openai/math/blob/main/preprints/Uniform-exclusion-of-Landau-Siegel-zeros-October-1-2026/paper.pdf}
}

@misc{Green2026Siegel,
  author       = {Green, Ben},
  title        = {On {Landau--Siegel} Zeros},
  year         = {2026},
  howpublished = {Hexagon preprint, hexagon:2610.00311v2},
  url          = {https://hexagonmath.org/2610.00311v2}
}

\end{document}